\documentclass[12pt,a4paper]{article}
\usepackage{amsmath}
\usepackage{amssymb}
\usepackage{amsthm}
\usepackage{amsfonts}
\usepackage{latexsym}

\theoremstyle{plain}
\newtheorem{thm}{Theorem}
\newtheorem{cor}{Corollary}

\theoremstyle{definition}

\theoremstyle{remark}

\title{A note on scaling asymptotics for Bohr-Sommerfeld Lagrangian submanifolds}
\author{Roberto Paoletti\footnote{\noindent{\bf Address:}
Dipartimento di Matematica e Applicazioni, Universit\`a degli Studi
di Milano Bicocca, Via R. Cozzi 53, 20125 Milano,
Italy; {\bf e-mail}: roberto.paoletti@unimib.it }}
\date{}

\begin{document}

\maketitle

\section{Introduction}

The purpose of this note is to improve an expansion  in \cite{dp}
for the asymptotics associated to Bohr-Sommerfeld Lagrangian submanifolds
of a compact Hodge manifold, in the context of geometric quantization
(see e.g. \cite{bw}, \cite{bpu}, \cite{gs-gc}, \cite{w}).
We adopt the general framework for quantizing Bohr-Sommerfeld Lagrangian submanifolds presented
in \cite{bpu}, based on applying the Szeg\"{o} kernel of the quantizing line bundle to certain
delta functions concentrated along the submanifold.

Let $M$ be a d-dimensional complex projective manifold, with complex structure $J$;
consider an ample line bundle $A$ on it, and let
$h$ be an Hermitian metric on $A$ such that the unique compatible connection has
curvature $\Omega=-2i\omega$, where $\omega$ is a K\"{a}hler form. Then the  unit
circle bundle $X\subseteq A^*$, endowed with the connection one-form $\alpha$, is a contact manifold.
A Bohr-Sommerfeld Lagrangian submanifold of $M$ (or, more precisely, of $(M,A,h)$)
is then simply a Legendrian submanifold $\Lambda\subseteq X$, conceived as an immersed submanifold of $M$.

In a standard manner, $X$ inherits a Riemannian structure for which the projection
$\pi:X\rightarrow M$ is a Riemannian fibration. In view of this, in the following at places we
shall implicitly identify (generalized) functions, densities and half-densities on $X$.

Referring to \S 2 of \cite{dp} for a more complete description of the prelimiaries involved,
we recall that if $\Lambda\subseteq X$ is a compact Legendrian submanifold, and $\lambda$ is a half-density on it,
there is a naturally induced generalized half-density $\delta_{\Lambda,\lambda}$ on $X$ supported on
$\Lambda$; following \cite{bpu}, we can then define a sequence of CR functions
$$
u_k=: \mathrm{P}_k\big(\delta_{\Lambda,\lambda}\big)\in \mathcal{H}(X)_k,
$$
where $\mathcal{H}(X)_k$ is the $k$-th isotype of the Hardy space with respect to the $S^1$-action,
and $\mathrm{P}_k:L^2(X)\rightarrow \mathcal{H}(X)_k$ is the orthogonal projector (extended to $\mathcal{D}'(X)\rightarrow
\mathcal{H}(X)_k$).
In the present setting there are natural unitary structures on $\mathcal{H}(X)_k$ and the space of global holomorphic
sections $H^0\left(M,A^{\otimes k}\right)$, and a natural unitary isomorphism
$\mathcal{H}(X)_k\cong H^0\left(M,A^{\otimes k}\right)$.
One thinks of the $u_k$'s as
representing the quantizations of $(\Lambda,\lambda)$ at Planck's constant $1/k$.
It is easily seen that $u_k$ is rapidly decaying as $k\rightarrow +\infty$ on
the complement of $S^1\cdot \Lambda=\pi^{-1}\Big(\pi(\Lambda)\Big)$; here $\pi:X\rightarrow M$ is the projection.
On the other hand, the asymptotic concentration of the $u_k$'s
along $S^1\cdot \Lambda$ poses an interesting problem, already considered in Theorem 3.12 of \cite{bpu}.

This theme was revisited in \cite{dp}, in a somewhat different technical setting; in particular,
Corollary 1.1 of \cite{dp} shows that the scaling asymptotics of $u_k$ (to be defined shortly)
near any $x\in S^1\cdot \Lambda$ admit an asymptotic expansion, and explicitly computes the leading order term.
We shall give presently a more precise description of this expansion, as a function on the tangent space of
$M$ at $m=\pi(x)$.
Namely, we shall show that this asymptotic expansion may be factored as an exponentially decaying term
in the component $w^\perp$ of $w\in T_mM$ orthogonal to $\Lambda$, times an asymptotic expansion with polynomial coefficients in $w$
(more precisely, the expansion is generally given by a finite sum of terms of this form, one from each branch of
$\Lambda$ projecting to $m$); the exponential term also
contains a symplectic pairing between $w^\perp$ and the component of $w$ along $\Lambda$, $w^\|$.
Furthermore, we shall give some relevant remainder estimates not mentioned in
\cite{dp}.

Before stating the results of this paper, let us recall that for any $x\in X$ we can find a
Heisenberg local chart for $X$ centered at $x$,
$$
\rho:B_{2\mathrm{d}}(\epsilon)\times (-\pi,\pi)\rightarrow X,\,\,\,\,\,\,(p,q,\theta)\mapsto
r_{e^{i\theta}}\Big(\mathfrak{e}\big(\varrho(p,q)\big)\Big);
$$
here $B_{2\mathrm{d}}(\epsilon)\subseteq \mathbb{R}^{2\mathrm{d}}$ is a ball of radius $\epsilon$ centered at the origin,
$\varrho:B_{2\mathrm{d}}(\epsilon)\rightarrow M$ is a preferred local chart for $M$ centered at $m=:\pi(x)$,
meaning that it trivializes the holomorphic and symplectic structures at $m$, and $\mathfrak{e}$ is a unitary local frame of
$A^*$, given by the unitarization of a preferred local frame
(complete definitions are in \cite{sz}).
Finally, $r_{e^{i\theta}}:X\rightarrow X$ is the diffeomorphism induced by the circle action.
It is in this kind of local coordinates that the scaling limits of Szeg\"{o} kernels exhibit their universal
nature (Theorem 3.1 of \cite{sz}). If $\rho$ is a system of Heisenberg local coordinates centered at $x$, and
$p,q\in \mathbb{R}^\mathrm{d}$, $w=p+iq$, one poses
$$
x+w=:\rho\big((p,q),0\big).
$$
For any $\theta$, we have
$$
u_k\big(\rho(p,q,\theta)\big)=e^{ik\theta}\,u_k\big(\rho(p,q,0)\big)=e^{ik\theta}\,u_k(x+w).
$$
Given that a system of Heisenberg local coordinates induces a unitary isomorphism
of $T_mM$ and $\mathbb{C}^\mathrm{d}$, with this understanding we can also consider the
expression $x+w$ with $w\in T_mM$.

If $x\in S^1\cdot \Lambda$, there are only finitely many elements $h_1,\ldots,h_{N_x}\in S^1$ such that
$x_j=:r_{h_j}(x)\in \Lambda$. Since $\Lambda$ is Legendrian, hence horizontal for the connection, for
any $j$ we may naturally identify the tangent space $T_{x_j}\Lambda\subseteq T_{x_j}X$ with a subspace
of $T_{\pi(x)}M$. With this in mind, if $w\in T_{\pi(x)}M$ we can write
$w=w_j^{\|}+w_j^\perp$ for unique $w_j^{\|}\in T_{x_j}\Lambda$ and $w_j^\perp\in T_{x_j}\Lambda^\perp$;
the latter denotes the orthocomplement of $T_{x_j}\Lambda$ in $T_{\pi(x)}M$ in the Riemannian metric of $M$.

Finally, let $\mathrm{dens}_\Lambda^{(1/2)}$ be the Riemannian half-density on $\Lambda$
(for the induced metric);
thus if $\lambda$ is a $\mathcal{C}^\infty$ half-density on $\Lambda$
we can write $\lambda=F_\lambda\cdot \mathrm{dens}^{(1/2)}_\Lambda$
for a unique $F_\lambda\in \mathcal{C}^\infty(\Lambda)$.

\begin{thm}
\label{thm:main-non-equiv}
Let $\Lambda\subseteq X$ be a compact Legendrian submanifold, and suppose
$\lambda$ is a smooth half-weight on it. For every $k=1,2,\ldots$, let $u_k=:\mathrm{P}_k\big(\delta_{\Lambda,\lambda}\big)$.
Suppose $x\in S^1\cdot \Lambda$,
and choose a system of Heisenberg local coordinates for $X$ù
centered at $x$. Let $h_1,\ldots,h_{N_x}\in S^1$ be the finitely many elements such that
$r_{h_j}(x)\in \Lambda$. Then:
\begin{enumerate}
\item Suppose $a>0$. Uniformly for $\min_j\big\{\|w_j^\perp\|\big\} \gtrsim k^{a}$, we have
  $$
  u_k\left(x+\frac{w}{\sqrt{k}}\right)=O\left (k^{-\infty}\right).
  $$

\item There exists polynomials $a_{lj}$ on $\mathbb{C}^\mathrm{d}$ such that
the following holds:
for $w\in T_{\pi(x)}M$ and $k,\ell=1,2,\ldots$, let us define
  \begin{eqnarray*}
  R_{k,\ell}\left(x,w\right)&=:&u_k\left(x+\frac{w}{\sqrt{k}}\right)\\
  &&-\left(\frac{2k}{\pi}\right)^{\mathrm{d}/2}
  \sum_{j=1}^{N_x}h_j^{-k}\,e^{-\|w_j^\perp\|^2-i\omega _{\pi(x)}(w_j^\perp,w_j^{\|})}\,F_\lambda(x_j)\cdot \left(1+\sum _{l=1}^\ell
  k^{-l/2}a_{lj}(w)\right).
  \end{eqnarray*}
  Then uniformly for $\|w\|\lesssim k^{1/6}$ we have
  \begin{equation}
  \label{eqn:estimate-remainder}
  \big |R_{k,\ell}\left(x,w\right)\big|\le C_\ell\,k^{(\mathrm{d}-\ell-1)/2}\,\sum_{j=1}^{N_x}e^{-\frac{1-\epsilon}{2}\|w_j^\perp\|^2}.
  \end{equation}
  \end{enumerate}
\end{thm}

\begin{cor} $\forall\,w\in T_{\pi(x)}M$, the following asymptotic expansion holds as $k\rightarrow +\infty$:
  $$
  u_k\left(x+\frac{w}{\sqrt{k}}\right)\sim \left(\frac{2k}{\pi}\right)^{\mathrm{d}/2}\,
  \sum_{j=1}^{N_x}h_j^{-k}\,e^{-\|w_j^\perp\|^2-i\omega _{\pi(x)}(w_j^\perp,w_j^{\|})}\,F_\lambda(x_j)\cdot \left(1+\sum _{l\ge 1}
  k^{-l/2}a_{lj}(w)\right).
  $$
  \end{cor}

This agrees with Corollary 1.1 of \cite{dp} to leading order, but gives a clearer picture of the asymptotic
expansion, as well as an explicit estimate on the remainder.

The proof of Theorem \ref{thm:main-non-equiv} is based on the scaling asymptotics of
Szeg\"{o} kernels in Theorem 3.1 of \cite{sz}, whereas the proofs in \cite{dp} are based on
microlocal arguments that encompass the equivariant setting, involving a direct use of the parametrix
developed in \cite{bs}. In view of the scaling asymptotics of
equivariant Szeg\"{o} kernels proved in \cite{p},
factorizations akin to Theorem \ref{thm:main-non-equiv}
also hold in the equivariant setting; we shall not discuss this here.

\section{Proof of Theorem \ref{thm:main-non-equiv}.}

Let us first prove 2. Thus, we want to investigate the asymptotics of $u_k\left(x+\frac{w}{\sqrt{k}}\right)$
as $k\rightarrow +\infty$, assuming that $w\in T_{\pi(x)}M$, $\|w\|\le C\,k^{1/6}$ for some fixed $C>0$.

Let $\Pi_k\in \mathcal{C}^\infty(X\times X)$ be the
Schartz kernel of $\mathrm{P}_k$; explicitly, if $\left\{s_r^{(k)}\right\}$
is an orthonormal basis of $\mathcal{H}(X)_k$, then
$$
\Pi_k(y,y')=\sum _rs_r^{(k)}(y)\cdot \overline{s_r^{(k)}(y')}\,\,\,\,\,\,\,\,\,\,\,(y,y'\in X).
$$
Let $\mathrm{dens}_X$ and $\mathrm{dens}_\Lambda$ denote, respectively, the Riemannian density
on $X$ and $\Lambda$.
Then, in standard distributional short-hand, by definition of $\delta_{\Lambda,\lambda}$
for any $x'\in X$
we have
\begin{eqnarray}
\label{eqn:integrale-su-lambda}
\lefteqn{u_k(x')=\int _X\Pi_k(x',y)\,\delta_{\Lambda,\lambda}(y)\,\mathrm{dens}_X(y)}\nonumber\\
&=&\left<\delta_{\Lambda,\lambda},
\Pi_k(x',\cdot)\right>=\int _\Lambda \Pi_k(x',y)\,F_\lambda(y)\,\mathrm{dens}_\Lambda(y).
\end{eqnarray}

Let us write $\mathrm{dist}_M$ for the Riemannian distance function on $M$, pulled-back to a smooth function
on $X\times X$ by the projection $\pi\times \pi$. Let us set:
\begin{eqnarray*}
V_k&=:&\left\{x'\in X\,:\,\mathrm{dist}_M\big(x,x')<4C\,k^{-1/3}\right\},\\
V_k'&=:&\left\{x'\in X\,:\,\mathrm{dist}_M\big(x,x')>3C\,k^{-1/3}\right\}.
\end{eqnarray*}
If $y\in V_k'$ and $\|w\|\le C\,k^{1/6}$, then
$\mathrm{dist}_M\left(x+\frac{w}{\sqrt{k}},y\right)\ge C\,k^{-1/3}$ for $k\gg 0$;
by the off-diagonal estimates on Szeg\"{o} kernels in \cite{christ}, therefore, $\Pi_k\left(x+\frac{w}{\sqrt{k}},y\right)=O\left(k^{-\infty}\right)$
uniformly for $y\in V_k'$.

For $k\gg 0$, $\Lambda\cap V_k$ has $N_x$ connected components:
$$
\Lambda\cap V_k=\bigcup _{j=1}^{N_x}\Lambda_{kj},
$$
where $\Lambda_{kj}$ is the connected component containing $x_j$.
Let $\{s_k,s'_k\}$ be an $S^1$-invariant partition of unity on $X$,
subordinate to the open cover
$\{V_k,V_k'\}$.
In view of (\ref{eqn:integrale-su-lambda}) and the previous discussion, we obtain
\begin{eqnarray}
\label{eqn:integrale-su-lambda-spezzato}
\lefteqn{u_k\left(x+\frac{w}{\sqrt{k}}\right)=\int _\Lambda \Pi_k\left(x+\frac{w}{\sqrt{k}},y\right)\,F_\lambda(y)\,\mathrm{dens}_\Lambda(y)}\nonumber\\
&\sim&\sum_{j=1}^{N_x}\int _{\Lambda_{kj}} \Pi_k\left(x+\frac{w}{\sqrt{k}},y\right)\,F_\lambda(y)\,s_{k}(y)\,\mathrm{dens}_\Lambda(y),
\end{eqnarray}
where $\sim$ means that the two terms have the same asymptotics.
Let us now evaluate the asymptotics of the $j$-th summand in (\ref{eqn:integrale-su-lambda-spezzato}).

To this end, recall that the Heisenberg local chart $\rho$ centered at $x$
depends on the choice of the preferred local chart
$\varrho$ at $\pi(x)$, and of the local frame $\mathfrak{e}$ of $A^*$.
We obtain a Heisenberg local chart $\rho'_j$ centered at $x_j$ by setting
$\rho'_j(p,q,\theta)=:r_{h_j}\Big(\rho(p,q,\theta)\Big)$.
By the discussion in \S 2 of \cite{dp} and (48) of the same paper,
we can compose $\rho'_j$ with a suitable transformation in $(p,q)$ (that is,
a change of preferred local chart for $M$) so as to obtain
a Heisenberg local chart $\rho_j(p,q,\theta)$ centered at $x_j$ with the following property:
$\Lambda$ is locally defined near $x_j$
by the conditions
$\theta=f_j(q)$ and $p=0$, where $f_j$ vanishes to third order at the origin.
By construction, we have
$\rho_j(p,q,\theta)=r_{h_j}\Big(\rho(p',q',\theta)\Big)$
for a certain local diffeomorphism $(p,q)\mapsto (p',q')$.

Thus $\Lambda$ is locally parametrized, near $x_j$ and in the chart $\rho_j$,
by the imaginary vectors $iq$;
viewing the $q$'s as local coordinates on $\Lambda$ near $x_j$, locally we have
$\mathrm{dens}_\Lambda=D_\lambda \cdot |dq|$, for a unique locally defined smooth function
$D_\Lambda$. By construction of Heisenberg local coordinates, $D_\lambda(0)=1$.

Applying a rescaling by $k^{-1/2}$, we obtain
\begin{eqnarray}
\label{eqn:integrale-su-lambda-spezzato-jth-term}
\lefteqn{\int _{\Lambda_{kj}} \Pi_k\left(x+\frac{w}{\sqrt{k}},y\right)\,F_\lambda(y)\,s_{k}(y)\,\mathrm{dens}_\Lambda(y)}\\
&=&k^{-\mathrm{d}/2}\,\int _{\mathbb{R}^\mathrm{d}}\,\Pi_k\left(x+\frac{w}{\sqrt{k}},
r_{e^{if_j\left(q/\sqrt{k}\right)}}\left(x_j+\frac{iq}{\sqrt{k}}\right)\right)\,
F_\lambda\left(\frac{q}{\sqrt{k}}\right)\,s_{k}\left(\frac{iq}{\sqrt{k}}\right)\,D_\lambda\left(\frac{q}{\sqrt{k}}\right)\, dq\nonumber\\
&=&k^{-\mathrm{d}/2}\,\int _{\mathbb{R}^\mathrm{d}}\,e^{-ikf_j\left(\frac{q}{\sqrt{k}}\right)}\,\Pi_k\left(x+\frac{w}{\sqrt{k}},
x_j+\frac{iq}{\sqrt{k}}\right)\,
F_\lambda\left(\frac{q}{\sqrt{k}}\right)\,s_{k}\left(\frac{iq}{\sqrt{k}}\right)\,D_\lambda\left(\frac{q}{\sqrt{k}}\right)\, dq.\nonumber
\end{eqnarray}
Here, $x+\frac{w}{\sqrt{k}}=\rho\left(\frac{\Re(w)}{\sqrt{k}},\frac{\Im(w)}{\sqrt{k}},0\right)$
(we use the Heisenberg chart to unitarily identify
$T_mM$ with $\mathbb{C}^\mathrm{d}$), and
$x_j+\frac{iq}{\sqrt{k}}=\rho_j\left(0,\frac{q}{\sqrt{k}},0\right)$.
Notice that
$s_{k}\left(\frac{iq}{\sqrt{k}}\right)=1$ for $\|q\|\lesssim k^{1/6}$,
$s_{k}\left(\frac{iq}{\sqrt{k}}\right)=0$ for $\|q\|\gtrsim k^{1/6}$.
In particular, integration takes place over a ball of radius
$\thicksim k^{1/6}$.
Also, Taylor expanding $F_\lambda$ and $f_j$ at the origin we have asymptotic expansions
$$
F_\lambda\left(\frac{q}{\sqrt{k}}\right)\sim
F_\lambda\left(x_j\right)+\sum _{r\ge 1} k^{-r/2}\,b_r(q),\,\,\,\,\,
D_\lambda\left(\frac{q}{\sqrt{k}}\right)\sim
1+\sum _{r\ge 1}k^{-r/2}\,c_r(q),
$$
and, since $f_j$ vanishes to third order at the origin,
$$
f_j\left(\frac{q}{\sqrt{k}}\right)\sim
\sum _{r\ge 0}k^{-(3+r)/2}\,d_r(q),\,\,\,\,\,
e^{-ikf_j\left(\frac{q}{\sqrt{k}}\right)}\sim 1+\sum _{r\ge 1}k^{-r/2}\,e_r(q),$$
for suitable polynomials $b_r$,$c_r$, $d_r$, $e_r$.

Let $w_j\in \mathbb{C}^\mathrm{d}$ correspond to $w$ in the Heisenberg local coordinates $\rho_j$.
By the above, Taylor expanding the transformation $(p,q)\mapsto (p',q')$, we
obtain $x+\frac{w}{\sqrt{k}}=r_{h_j^{-1}}\left (x_j+\frac{w_j}{\sqrt{k}}+H(w,k)\right)$,
where $H(w,k)\sim \sum _{f\ge 2}k^{-f/2}h_f(w)$.
Without affecting the leading order term of the resulting asymptotic expansion, we may pretend for simplicity that
$x+\frac{w}{\sqrt{k}}=r_{h_j^{-1}}\left (x_j+\frac{w_j}{\sqrt{k}}\right)$.

Write $w_j=p_j+iq_j$, with $p_j,q_j\in \mathbb{R}^\mathrm{d}$. Thus $w_j^{\perp}=p_j$,
$w_j^\|=iq_j$.
In view of Theorem 3.1 of \cite{sz}, we have
\begin{eqnarray}
\label{eqn:asympt-exp-for-pi-equiv}
\lefteqn{\Pi_k\left(x+\frac{w}{\sqrt{k}},
x_j+\frac{iq}{\sqrt{k}}\right)}\\
&=&\Pi_k\left(r_{h_j^{-1}}\left (x_j+\frac{w_j}{\sqrt{k}}\right),
x_j+\frac{iq}{\sqrt{k}}\right)=h_j^{-k}\,\Pi_k\left(x_j+\frac{w_j}{\sqrt{k}},
x_j+\frac{iq}{\sqrt{k}}\right)\nonumber\\
&\sim&h_j^{-k}\,\left(\frac{k}{\pi}\right)^{\mathrm{d}}\,e^{-ip_j\cdot q-\frac 12\|p_j\|^2
-\frac 12\|q_j-q\|^2}\cdot \left(1+
\sum_{r\ge 1}k^{-r/2}R_j(w,q)\right),\nonumber
\end{eqnarray}
for certain polynomials $R_j$ in $w$ and $q$.
Furthermore, by the large ball estimate on the remainder discussed in \S 5 of \cite{sz},
the remainder after summing over $1\le r\le R$ is bounded by
\begin{equation}
\label{eqn:remainder-sz}
C_R\,k^{\mathrm{d}-(R+1)}\,e^{-\frac{1-\epsilon}{2}\,(\|p_j\|^2+\|q-q_j\|^2)}.
\end{equation}
It follows that the product of these asymptotic expansions can be integrated term by term; given this,
we only lose a contribution which is $O\left(k^{-\infty}\right)$ by setting $s_k=1$ and
integrating over all of $\mathbb{R}^\mathrm{d}$.

We have
$$
\int _{\mathbb{R}^\mathrm{d}}e^{-ip_j\cdot q
-\frac 12\|q_j-q\|^2}\,dq=e^{-ip_j\cdot q_j}\int _{\mathbb{R}^\mathrm{d}}e^{-ip_j\cdot s-\frac 12\|s\|^2}\,ds=
(2\pi)^{\mathrm{d}/2}\,e^{-ip_j\cdot q_j-\frac 12\|p_j\|^2}.
$$
Given (\ref{eqn:asympt-exp-for-pi-equiv}), this implies that (\ref{eqn:integrale-su-lambda-spezzato-jth-term})
is given by an asymptotic expansion, with leading order term
$$
h_j^{-k}\,\left(\frac{2k}{\pi}\right)^{\mathrm{d}/2}\,e^{-\|w_j^\perp\|^2-i\omega _{\pi(x)}(w_j^\perp,w_j^{\|})}
\,F_\lambda\left(x_j\right).
$$
To determine the general term of the expansion, on the other hand, we are led to computing integrals of the form
$$
\int _{\mathbb{R}^\mathrm{d}}q^\beta e^{-ip_j\cdot q
-\frac 12\|q_j-q\|^2}\,dq=e^{-ip_j\cdot q_j}\int _{\mathbb{R}^\mathrm{d}}(s+q_j)^\beta\,
e^{-ip_j\cdot s-\frac 12\|s\|^2}\,ds.
$$
where $q^\beta$ is some monomial.
Thus we led to a sum of terms of the form
$$
e^{-ip_j\cdot q_j}\,C_\gamma(q_j)\,\int _{\mathbb{R}^\mathrm{d}}s^\gamma\,
e^{-ip_j\cdot s-\frac 12\|s\|^2}\,ds,
$$
and the integral is the evaluation at $p_j$ of the Fourier transform of $s^\gamma\,
e^{-\frac 12\|s\|^2}$. Up to a scalar factor, the latter is an iterated derivative to
$e^{-\frac 12\|s\|^2}$; therefore we are left with a summand of the form
$
e^{-ip_j\cdot q_j}\,T(q_j,q_j)\,
e^{-\frac 12\|p_j\|^2}\,ds$,
where $T$ is a polynomial in $p_j,q_j$.
Given (\ref{eqn:asympt-exp-for-pi-equiv}), this implies that the general term of the asymptotic expansion
for (\ref{eqn:integrale-su-lambda-spezzato-jth-term})
has the form
$$
h_j^{-k}\,\left(\frac{2k}{\pi}\right)^{\mathrm{d}/2}\,k^{-l/2}\,e^{-\|w_j^\perp\|^2-i\omega _{\pi(x)}(w_j^\perp,w_j^{\|})}
\,F_\lambda\left(x_j\right)\cdot a_{lj}(w)
$$
for an appropriate polynomial $a_{lj}(w)$.
Finally, (\ref{eqn:estimate-remainder}) (at $x_j$) follows by integrating (\ref{eqn:remainder-sz}).

To complete the proof of 2., we need only sum over $j$.

Turning to the proof of 1., by definition of preferred local coordinates, if
$w^{\perp}_j\ge C\,k^a$, say, then
$$
\mathrm{dist}_M\left (x+\frac{v}{\sqrt{k}},\Lambda_{kj}\right)\ge \frac{C}{2}\,k^{a-\frac 12},
$$
for all $k\gg 0$.
By the off-diagonal estimates of \cite{christ},
$\Pi_k\left(x+\frac{w}{\sqrt{k}},y\right)=O\left(k^{-\infty}\right)$ uniformly for
$y\in \Lambda_{kj}$.

\hfill Q.E.D.

\end{document}